\documentclass[10pt]{amsart}

\usepackage{amsmath,amssymb,amsthm,a4wide}

\usepackage{tikz}
\usetikzlibrary{shapes}
\usetikzlibrary{intersections}

\usepackage{graphicx}
\usepackage{graphics}

\newtheorem*{thm}{Theorem}
\newtheorem{lemma}{Lemma}

\newcommand{\supp}{\operatorname{supp}}

\begin{document}

\title[]{On Suprema of Autoconvolutions \\with an application to sidon sets}
\author{Alexander Cloninger}
\address[Alexander Cloninger]{Department of Mathematics, Program in Applied Mathematics, Yale University, New Haven, CT 06510, USA}
\email{alexander.cloninger@yale.edu}

\author{Stefan Steinerberger}
\address[Stefan Steinerberger]{Department of Mathematics, Yale University, 06511 New Haven, CT, USA}
\email{stefan.steinerberger@yale.edu}

\begin{abstract} 
Let $f$ be a nonnegative function supported on $(-1/4, 1/4)$. We show
$$ \sup_{x \in \mathbb{R}}{\int_{\mathbb{R}}{f(t)f(x-t)dt}} \geq 1.28\left(\int_{-1/4}^{1/4}{f(x)dx} \right)^2,$$
where 1.28 improves on a series of earlier results. 
The inequality arises naturally in additive 
combinatorics in the study of Sidon sets. We derive a relaxation of the problem that reduces to a finite number of cases and yields slightly stronger results. Our approach should be able to prove lower bounds that are arbitrary close to the sharp result. Currently, the bottleneck in our approach is runtime: new ideas might be able to significantly speed up the computation.
\end{abstract}
\maketitle

\section{Introduction}
\subsection{Sidon sets.}
A subset $A \subset \left\{1, 2, \dots, N\right\}$ is called $g-$Sidon if 
$$ | \left\{(a,b) \in A\times A: a+b=m\right\}| \leq g$$
for every $m$. How large can $g-$Sidon sets possibly be? Denote the answer by
$$ \beta_g(n) := \max_{A \subset \left\{1, 2, \dots, N\right\} \atop \mbox{\tiny A is g-Sidon}}{|A|}.$$
There has been a lot of research activity on bounding these quantities: we only deal with the asymptotic
case as $g$ becomes large. Cilleruelo, Ruzsa \& Vinuesa \cite{cill} have recently shown that the maximal cardinality of a $g-$Sidon set satisfies
$$\underline{\sigma(g)}\sqrt{gn}(1-o(1)) \leq \beta_g(n) \leq \overline{\sigma(g)}\sqrt{gn}(1+o(1)),$$
where the $o(1)$ is with respect to $n$ and 
$$ \lim_{g \rightarrow \infty}{\underline{\sigma(g)}} = \sigma = \lim_{g \rightarrow \infty}{\overline{\sigma(g)}}$$
for some universal constant $\sigma \in \mathbb{R}$. 

\subsection{An equivalent continuous problem.} As proven by Cilleruelo, Ruzsa \& Vinuesa \cite{cill}, the constant $\sigma$ has an alternative representation as the solution of a continuous problem involving the autoconvolution $f*f$ which was first considered by Schinzel \& Schmidt \cite{sch}. Throughout this paper, we will use
$$ (f*g)(x) = \int_{\mathbb{R}}{f(t)g(x-t)dt}$$
to denote the convolution. Consider all nonnegative functions $f$ supported on $[0,1]$ satisfying
$$\|f*f\|_{L^{\infty}(\mathbb{R})} =  \sup_{x \in \mathbb{R}}{\int_{\mathbb{R}}^{}{f(t)f(x-t)dt}} \leq 1,$$
then we have
$$ \int_{\mathbb{R}}{f(x)dx} \leq \sigma \qquad \mbox{and this inequality is sharp.}$$
A rephrasing of the statement is as follows: let $f$ be a nonnegative
 function supported on the interval $[-1/4, 1/4]$. Then we have that
$$ \sup_{x \in \mathbb{R}}{\int_{\mathbb{R}}{f(t)f(x-t)dt}} \geq \frac{2}{\sigma^2}\left(\int_{-1/4}^{1/4}{f(x)dx} \right)^2.$$
We abbreviate henceforth $c = 2/\sigma^2.$ As for the lower bound,
we have
\begin{align*}
c &\geq 1\qquad &&\mbox{trivial}\\
&\geq 1.151 \qquad &&\mbox{Cilleruelo, Ruzsa \& Trujillo \cite{cill1}} \\
&\geq  1.178 &&\mbox{Green \cite{green}}\\
&\geq  1.183 &&\mbox{Martin \& O'Bryant \cite{martin}}\\
&\geq  1.251 &&\mbox{Yu \cite{yu}}\\
&\geq  1.263 &&\mbox{Martin \& O'Bryant \cite{martin2}}\\
&\geq  1.274 &&\mbox{Matolcsi \& Vinuesa \cite{mat}}
\end{align*}
Work on the upper bound has been
carried out by Cilleruelo, Ruzsa \& Trujillo \cite{cill1}, Kolountzakis \cite{kol}, Martin \& O'Bryant \cite{martin} and Cilleruelo \& Vinuesa \cite{cill2}. 
The current state of the art is
$$ 1.2749 \leq c \leq 1.50992$$
with both bounds coming from Matolcsi \& Vinuesa \cite{mat}. The upper bound 1.50992 is conjectured to be almost thight and disproves the conjecture $c = \pi/2$ 
due to Schinzel \& Schmidt \cite{sch}. The upper bound comes from an explicit example; the lower 
bound is established using Fourier methods and earlier arguments of 
Martin \& O'Bryant \cite{martin3} and Yu \cite{yu}. Matolcsi \& Vinuesa also claim the theoretical limit of their approach for
 the lower bound to be at 1.276. 

\subsection{The result.} The purpose of this paper is to improve on the lower bound by proving that there is a relaxed problem which can be dealt with computationally. 
The actual improvement is small (though comparable to previous improvements). We believe the true merit of our contribution to be in demonstrating a new approach to the problem that might prove
quite effective.

\begin{thm} Let $f:\mathbb{R} \rightarrow \mathbb{R}_{\geq 0}$ be supported in $[-1/4,1/4]$. Then
$$ \sup_{x \in \mathbb{R}}{\int_{\mathbb{R}}{f(t)f(x-t)dt}} \geq 1.28 \left(\int_{-1/4}^{1/4}{f(x)dx} \right)^2.$$
\end{thm}
Our method seems to only be limited by our ability to do large-scale computations and we consider it likely that variants of our idea might be
much faster and able to drastically improve on the result outlined here.
There are two main ingredients: the first is to regard convolution
as a bilinear operator that respects spatial decompositions via
$$ \supp(f*g) \subseteq \supp(f) + \supp(g)$$
and acts nicely with respect to the $L^1-$norm via
$$ \int_{\mathbb{R}}{(f*g)(x)dx} = \left( \int_{\mathbb{R}}{f(x)dx}\right) \left( \int_{\mathbb{R}}{g(x)dx}\right)$$
for sufficiently regular functions.
The second
ingredient is the pigeonhole principle stating that for any $h:I \rightarrow \mathbb{R}_{}$
$$ \|h\|_{L^{\infty}} \geq \frac{1}{|I|}\int_{I}{|h(x)| dx}.$$

\section{Some intuition} 
There is a trivial proof of $c \geq 1$: since $f$ is supported
on $(-1/4, 1/4)$, the function $f*f$ must be supported on $(-1/2, 1/2)$. Furthermore, by Fubini, we have that
$$\int_{\mathbb{R}}^{}{(f*f)(x)dx} = \int_{\mathbb{R}}{\int_{\mathbb{R}}^{}{f(t)f(x-t)dt}dx} =
\left(\int_{\mathbb{R}}^{}{f(x)dx}\right)^2$$
and thus, trivially
$$ \left(\int_{\mathbb{R}}^{}{f(x)dx}\right)^2 = \|f*f\|_{L^1(-1/2,1/2)} \leq \|f*f\|_{L^{\infty}(-1/2,1/2)}.$$
This argument inspired our approach: instead of trying to show that $\|f*f\|_{L^{\infty}(\mathbb{R})}$ has
to be large by showing that $f *f$ has to be large in a single point, we show that there has to exist an interval
$I \subset (-1/2, 1/2)$ such that $\|f*f\|_{L^1(I)}$ is large and then use
$$ \frac{1}{|I|}\|f*f\|_{L^1(I)} \leq \|f*f\|_{L^{\infty}(I)} \leq \|f*f\|_{L^{\infty}(\mathbb{R})}.$$
Indeed, this also implies that we will actually be able to deduce slightly stronger results. Our main
result reads as follows and clearly implies the bound $c \geq 1.28$.
\begin{thm}[Full result] Let $f:\mathbb{R} \rightarrow \mathbb{R}_{\geq 0}$ be supported in $[-1/4,1/4]$. Then there
exists an interval $J \subset (-1/2, 1/2)$ of length $|J| \geq 1/48$ such that
$$ \int_{J}{(f*f)(x)dx} \geq 1.28 |J|  \left(\int_{-1/4}^{1/4}{f(x)dx} \right)^2.$$
\end{thm}

We will now describe the argument in greater detail:
let us decompose the interval $(-1/4,1/4)$ into four intervals $I_1, I_2, I_3, I_4$ of equal size in the
canonical way, i.e.
$$ I_j = \left(\frac{-3+j}{8}, \frac{-2+j}{8}\right) \qquad \mbox{for} \quad j = 1,2,3,4.$$
 Denoting the characteristic function of an interval $I$ by $\chi_I$, we define
$$f_j = f\chi_{I_j}$$
and use the linearity of the convolution to write
\begin{align*}
 f*f &= \left(\sum_{j=1}^{4}{f_j}\right)* \left(\sum_{j=1}^{4}{f_j}\right) \\
&= \sum_{j,k = 1}^{4}{f_j * f_k}
\end{align*}

\begin{figure}[h!]
\begin{tikzpicture}[xscale = 0.7, yscale = 0.7]
\draw [thick] (0,0) -- (10,0);
\draw[very thick] (0,0) to [out=00,in=195] (2,1.5);
\node at (1.3, -0.4) {$I_1$};
\draw[very thick] (2,1.5) to [out=15,in=190] (5,3);
\node at (3.5, -0.4) {$I_2$};
\draw[very thick] (5,3) to [out=10,in=170] (7,2);
\node at (6.3, -0.4) {$I_3$};
\draw[very thick] (7,2) to [out=-10,in=180] (10,0);
\node at (8.5, -0.4) {$I_4$};
\draw (-0.3,-0.4) node{$-\frac{1}{4}$};
\draw [fill] (2.5,0) circle [radius=0.05];
\draw (2.5,-0.4) node{$-\frac{1}{8}$};
\draw [fill] (0,0) circle [radius=0.05];
\draw (10.3,-0.4) node{$\frac{1}{4}$};
\draw [fill] (5,0) circle [radius=0.05];
\draw [fill] (10,0) circle [radius=0.05];
\draw (5,-0.4) node{$0$};
\draw [fill] (7.5,0) circle [radius=0.05];
\draw (7.5,-0.4) node{$\frac{1}{8}$};
\draw (2.5,0) -- (2.5,1.7);
\draw (5,0) -- (5,3);
\draw (7.5,0) -- (7.5,1.8);
\end{tikzpicture}
\caption{A nonnegative function supported on $(-1/4,1/4)$.}
\end{figure}
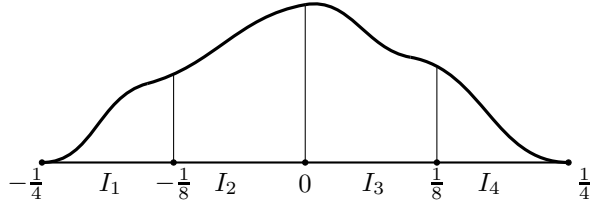

It is easy to see that $f_j * f_k$ is supported on $I_{j} + I_{k}$, where $+$ is interpreted as the
Minkowski sum of sets. This property allows now to deduce precisely the location of the support of each
term. We have, additionally, by Fubini, that
$$ \int_{\mathbb{R}}{(f_j * f_k)(x)dx} = \left(\int_{\mathbb{R}}{f_j(x)dx}\right)\left(\int_{\mathbb{R}}{f_k(x)dx}\right).$$
This allows now to deduce, for example, the pointwise inequality 
$$  (f_1*f_1)(x) + (f_1*f_2)(x) + (f_2*f_1)(x) \leq (f*f)(x) \qquad \mbox{on the interval} \quad (-1/2,-1/8).$$
Integrating the inequality on both sides yields
$$ \left(\int_{\mathbb{R}}{f_1(x)dx}\right)^2 + 2\left(\int_{\mathbb{R}}{f_1(x)dx}\right)\left(\int_{\mathbb{R}}{f_2(x)dx}\right) 
\leq \int_{-\frac{1}{2}}^{-\frac{1}{8}}{(f*f)(x)dx} \leq \frac{3}{8}\|f*f\|_{L^{\infty}(\mathbb{R})}.$$
This corresponds to certain restrictions on the distribution of the $L^1-$mass of the function
over the intervals assuming $\|f*f\|_{L^{\infty}(\mathbb{R})}$ to be small, or, arguing conversely,
shows that any function with $\|f*f\|_{L^{\infty}(\mathbb{R})}$ small must induce a partition of
its $L^1-$mass that satisfies the particular inequality. The idea behind our approach is mainly to check
whether it is at all possible for a function to satisfy a long list of such inequalities.
Note that the approach
is quite insensitive as to what the function actually looks like -- the only information used is the integral
over a small interval: this implies that it suffices to check the chain of inequalities for all step functions.

\section{Details and Proofs}
\subsection{Step functions.} This section is devoted to a simple Lemma capturing the main idea on a single length scale. It is 
part of the argument that we do not use any information on the function excepts its integral over
certain intervals -- the actual shape of the function over that interval plays no role. This has the effect
 that we may later restrict ourselves to step functions. 
\begin{lemma} For any $n \in \mathbb{N}$ let
$$ A_n = \left\{(a_{-n}, a_{-n+1}, \dots, a_{n-1}) \in (\mathbb{R}_{+})^{2n}: \sum_{i=-n}^{n-1}{a_i} = 4n\right\}$$
and
$$a_n := \min_{a \in A_n} \max_{2 \leq \ell \leq 2n} ~~~~\max_{-n \leq k \leq n-\ell}~~
\frac{1}{4n\ell}\sum_{k \leq i+j \leq k +\ell - 2}{a_i a_j},$$
where $k, l \in \mathbb{Z}$. Then
$$ c \geq a_n.$$
\end{lemma}

\begin{proof} Let $\varepsilon > 0$ be arbitrary and let $f$ be a 
nonnegative function supported on the interval $(-1/4, 1/4)$ with normalized mean 
$$\int_{-1/4}^{1/4}{f(x)dx} = 1$$
such that
$$ \sup_{x \in \mathbb{R}}{\int_{\mathbb{R}}{f(t)f(x-t)dt}} \leq c + \varepsilon.$$
Consider the decomposition of the interval into $2n$ equally
sized intervals
$$ \left( -\frac{1}{4}, \frac{1}{4}\right) = \bigcup_{j=-n}^{n-1}{I_j} \qquad \mbox{where}\quad
I_j = \left(\frac{j}{4n},\frac{j+1}{4n}\right).$$
We denote the restriction of $f$ to the interval $I_j$ by $f_j$ and define the average on that region by
$$ a_{j} = \frac{1}{|I_j|}\int_{I_j}{f(x)dx} = 4n\int_{\mathbb{R}}{f_j(x)dx}.$$
Trivially, since $f \geq 0$, we have $a_{j} \geq 0$ as well as
$$ \frac{1}{4n}\sum_{j=-n}^{n-1}{a_j} = \int_{-1/4}^{1/4}{f(x)dx} = 1$$
and thus $(a_{-n}, a_{-n+1}, \dots, a_{n-1}) \in A_n$. Let $2 \leq \ell \leq 4n$ be an arbitrary integer and let $-2n \leq k \leq 2n-\ell$ be another arbitrary integer. A simple computation yields

\begin{align*}
\frac{1}{4n\ell}\sum_{k \leq i+j \leq k +\ell -2}{a_i a_j} &=
\frac{4n}{\ell}\sum_{k \leq i+j \leq k +\ell -2}{\frac{a_i}{4n} \frac{a_j}{4n}} = \frac{4n}{\ell}\sum_{k \leq i+j \leq k +\ell -2}{\left(\int_{I_i}{f(x)dx}\right)\left( \int_{I_j}{f(x)dx}\right)} \\
&=\frac{4n}{\ell}\sum_{k \leq i+j \leq k +\ell -2}{\int_{\mathbb{R}}{(f_i*f_j)(x)dx}} \underbrace{\leq}_{(\diamond)} \frac{4n}{\ell}\int_{\frac{k}{4n}}^{\frac{k+l}{4n}}{(f*f)(x)dx} \\
&\leq \frac{4n}{\ell}\int_{\frac{k}{4n}}^{\frac{k+l}{4n}}{\|f*f\|_{L^{\infty}(\mathbb{R})}dx} = \|f*f\|_{L^{\infty}} \leq c+\varepsilon,
\end{align*}
where $(\diamond)$ is the only nonobvious step, which follows from the fact that $k \leq i+j \leq k +\ell -2$
implies that
$$ I_{i} + I_{k} \subseteq \left( \frac{k}{4n}, \frac{k+l}{4n} \right).$$
Since $\varepsilon$ was arbitrary, the compactness of $A_n$ now implies
$$\min_{a \in A_n} ~~~~\max_{2 \leq \ell \leq 2n} ~~~~\max_{-n \leq k \leq n-\ell}~~
\frac{4n}{\ell}\sum_{k \leq i+j \leq k +\ell}{a_i a_j} \leq c.$$
\end{proof}

\subsection{Approximation of the simplex.} We will now see that the simplex $A_n$ can be suitably discretized. This problem has been
considered at a greater level of generality in the literature \cite{bonze}, however, we only require a special case
and can thus provide a simpler self-contained argument.

\begin{lemma} For any $n, m \in \mathbb{N}$ let
\begin{align*}
A_n &= \left\{(a_{-n}, a_{-n+1}, \dots, a_{n-1}) \in (\mathbb{R}_{+})^{2n}: \sum_{i=-n}^{n-1}{a_i} = 4n\right\}\\
B_{n,m} &= \left\{(b_{-n}, b_{-n+1}, \dots, b_{n-1}) \in \left(\frac{1}{m}\mathbb{N}\right)^{2n}: \sum_{i=-n}^{n-1}{b_i} = 4n\right\}.
\end{align*}
Then $B_{n,m}$ is a $1/m^{}-$net of $A_{n}$ in the $\ell^{\infty}-$norm, i.e.
$$ \forall a \in A_n~~\exists ~b \in B_{n,m} \qquad \quad \| a - b\|_{\ell^{\infty}} \leq \frac{1}{m}.$$
\end{lemma}
\begin{proof} We give an explicit construction rule for $b$. First, for all
$-n \leq i \leq n-1$ with 
$$a_i \in \frac{1}{m}\mathbb{N} \quad \mbox{we set} \quad b_i = a_i.$$
Since $a_i \geq 0$, this implies that all entries \textit{not} of that form are positive and therefore we can
either round up or round down and choose either
$$     b_i = \frac{\lfloor m a_i \rfloor}{m}    \qquad \mbox{or} \qquad    b_i = \frac{\lceil m a_i \rceil}{m}$$
without violating the desired bound $|a_i - b_i| \leq m^{-1}$.
We now need to define an explicit rule that tells us whether to round up or round down in every single
instance. This choice is decided as follows: the first time we have to decide, we round down. Note that if we have to
decide once, we have to decide at least twice. We define a variable
$$ t_{i} = \sum_{j < i}{(b_j - a_j)}.$$
Since we round down in the first instance, the first nonzero value in the variable $t$ is going to be negative but bigger
than $-1/m$.
We will then, in each instance $-n \leq i \leq n-1$ where a decision about rounding up or down is required,
proceed in a way that ensures $-1/m < t_{i+1} \leq 0$. This is always possible: if rounding up $a_i$ were to lead
to 
$$ t_{i+1} > 0,$$
then, by definition,
$$ 0 < t_{i+1} = t_{i} + \frac{\lceil m a_i \rceil}{m} - a_i =  t_{i} +\frac{\lfloor m a_i \rfloor}{m} + \frac{1}{m} - a_i$$
in which case, we see that rounding down leads to a value
$$ t_{i+1} =  t_{i} + \left( \frac{\lfloor m a_i \rfloor}{m} - a_i \right) > - \frac{1}{m}.$$
This construction rule now proceeds a series of nonnegative rationals that satisfy $|a_i - b_i| \leq m^{-1}$. It
remains to see why the rule ends up selecting a series of rationals that add up to the same number $4n$.
Clearly, this requires us to set
$$ b_{n-1} := 4n - \sum_{j=-n}^{n-2}{b_j}$$
but it is not a priori clear whether this guarantees $b_{n-1} \geq 0$ or $|b_{n-1} - a_{n-1}| \leq m^{-1}$. Multiplying the definition of $b_{n-1}$ by $m$ implies immediately that $m b_{n-1} \in \mathbb{N}$
and we will now show the two remaining properties to be satisfied.
First we remark that, again by construction,
$$ - \frac{1}{m} \leq  \sum_{j =1}^{n-1}{(b_j - a_j)} \leq 0$$
which implies
$$ b_{n-1} = 4n - \sum_{j=-n}^{n-2}{b_j} = a_{n-1} + \sum_{j=-n}^{n-2}{a_j}  -  \sum_{j=-n}^{n-2}{b_j}  = a_{n-1} - t_{n-1}.$$
This yields $b_{n-1} \geq a_{n-1} \geq 0$ as well as $|b_{n-1} - a_{n-1}| \leq |t_{n-1}| < m^{-1}$.
\end{proof}

\begin{lemma}[Discretization] \label{disc}  For any $n, m \in \mathbb{N}$ let
\begin{align*}
B_{n,m} &= \left\{(b_{-n}, b_{-n+1}, \dots, b_{n-1}) \in \left(\frac{1}{m}\mathbb{N}\right)^{2n}: \sum_{i=-n}^{n-1}{b_i} = 4n\right\} 
\end{align*}
and
\begin{align*}
b_{n,m} &:= \min_{b \in B_{n,m}} \max_{2 \leq \ell \leq 2n} ~~~~\max_{-n \leq k \leq n-\ell}~~
\frac{1}{4n\ell}\sum_{k \leq i+j \leq k +\ell - 2}{b_i b_j}.
\end{align*}
 Then
$$  c \geq b_{n,m} - \frac{2}{m} - \frac{1}{m^2}.$$
\end{lemma}
\begin{proof} The first Lemma allows us to replace a general function by a vector $a \in A_{n}$, which - in turn - we interpret as a step function $f$. The second Lemma guarantees
the existence of $b \in B_{n,m}$ approximating $a$ which we interpret as another step-function $g$
satisfying
$$ \|f-g\|_{L^{\infty}(\mathbb{R})} \leq \frac{1}{m}.$$
Alternatively, we can write that $f(x) = g(x) + \varepsilon(x),$ where $\varepsilon(x)$ is a third step function that is uniformly bounded by $m^{-1}$. 
By associativity of the convolution, we have that
$$ (g*g) = (f*f) - 2(f*\varepsilon) + (\varepsilon*\varepsilon).$$
We have 
$$ | (\varepsilon*\varepsilon)(x)|= \left|\int_{-\frac{1}{4}}^{\frac{1}{4}}{\varepsilon(x-y)\varepsilon(y) dy} \right| \leq \frac{1}{m}\int_{-\frac{1}{4}}^{\frac{1}{4}}{ |\varepsilon(x)| dx} \leq  \frac{1}{m^2}.$$
As for the other term, we have the inequality
$$ |(f*\varepsilon)(x)| = \left|\int_{-\frac{1}{4}}^{\frac{1}{4}}{f(x-y)\varepsilon(y) dy} \right| \leq  \frac{1}{m}\int_{-\frac{1}{4}}^{\frac{1}{4}}{f(y) dy} \leq \frac{1}{m}.$$
Using Lemma 1, this proves the result.
\end{proof}

\subsection{Conclusion of the argument} 
The ideas above were additionally coupled with a multi-scale argument that has no additional mathematical content but speeds up the implementation. The guiding idea is that for some distributions of $L^1-$mass on $(-1/4,1/4)$ it is sufficient to look at very rough discretizations
whereas other distributions will require a finer analysis. Consider, for example, a function $f$ as depicted in Fig. 2. Assuming it to be $L^1-$normalized, we can see that at least $80\%$ of the $L^1-$mass are contained in $(-1/4, 0)$ and therefore at
least $64\%$ of the $L^1-$mass of $(f*f)$ is contained in $(-1/2, 0)$ which already implies
$$ \sup_{-1/2 \leq x \leq 0}{(f*f)(x)} \geq \frac{0.64}{0.5} = 1.28.$$
If our goal is to prove a lower bound of $1.28$, we can effectively restrict our search for counterexamples to functions satisfying
$$ 0.2 \leq \int_{-\frac{1}{4}}^{0}{f(x)dx} \leq 0.8 \qquad \mbox{and, by symmetry,} \qquad 0.2 \leq \int_{0}^{\frac{1}{4}}{f(x)dx} \leq 0.8.$$
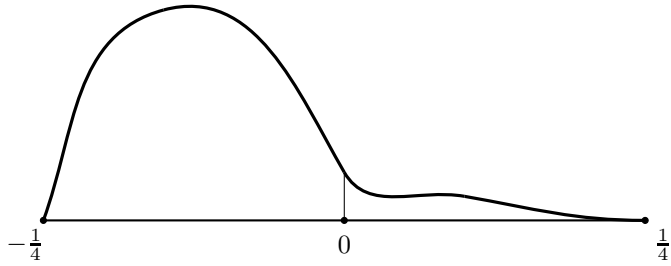
\begin{figure}[h!]
\begin{tikzpicture}[xscale = 0.8, yscale = 0.8]
\draw [thick] (0,0) -- (10,0);
\draw[very thick] (0,0) to [out=70,in=195] (2,3.5);
\draw[very thick] (2,3.5) to [out=15,in=120] (5,0.8);
\draw[very thick] (5,0.8) to [out=300,in=170] (7,0.4);
\draw[very thick] (7,0.4) to [out=-10,in=180] (10,0);
\draw (-0.3,-0.4) node{$-\frac{1}{4}$};
\draw [fill] (0,0) circle [radius=0.05];
\draw (10.3,-0.4) node{$\frac{1}{4}$};
\draw [fill] (5,0) circle [radius=0.05];
\draw [fill] (10,0) circle [radius=0.05];
\draw (5,-0.4) node{$0$};
\draw (5,0) -- (5,0.8);
\end{tikzpicture}
\caption{A nonnegative function supported on $(-1/4,1/4)$.}
\end{figure}

Suppose now there exists a function with $\|f * f \|_{L^{\infty}} \leq 1.28$. Then Lemma \ref{disc} implies that
$$ \forall n \in \mathbb{N} \qquad b_{n,m} \leq 1.28 + \frac{2}{m} + \frac{1}{m^2}$$
and we will now disprove this.  While implementing this scheme, we realized that it is quite a bit more effective in Lemma 3 to not use
$$ |(f*\varepsilon)(x)| = \left|\int_{-\frac{1}{4}}^{\frac{1}{4}}{f(x-y)\varepsilon(y) dy} \right| \leq  \frac{1}{m}\int_{-\frac{1}{2}}^{\frac{1}{2}}{f(x-y) dy} \leq \frac{1}{m}$$
but a refined version exploiting the support more effectively via
\begin{eqnarray}\label{eq:refinedLowerBound}
|(f*\varepsilon)(x)| = \left|\int_{-\frac{1}{4}}^{\frac{1}{4}}{f(x-y)\varepsilon(y) dy} \right|  \leq  \frac{1}{m}\int_{-\frac{1}{4}+\max(x,0)}^{\frac{1}{4}+\min(0,x)}{f(x-y) dy}
\end{eqnarray}
which is more efficient at excluding cases. We used $m=50$ and started with $n=3$. The same argument as above allows to ensure that some configurations in $B_{3,50}$ cannot come from
a function satisfying $\|f * f \|_{L^{\infty}} \leq 1.28$ but other configurations cannot be ruled out at such a rough scale. All cases that cannot be ruled out that level are being refined
by a dyadic split of each interval. This generates elements in $B_{6,50}$ where again some can be ruled out to come from a function $\|f * f \|_{L^{\infty}} \leq 1.28$ while others
cannot. We repeat the scheme as long as required. The final outcome was that no element in $B_{24,50}$ can come from a function $f$ satisfying $\|f * f \|_{L^{\infty}} \leq 1.28$ which implies that such a function cannot exist.  $\qed$

\subsection{Implementation}
We decided to run these cases in parallel, and moreover to utilize graphical processing units (GPUs).  GPU computing greatly reduce the cost of computation, as GPUs are highly optimized for vector operations.  However, there is overhead cost to transferring data to GPUs, and GPUs are suboptimal for any computation other than matrix multiplication. This motivated us to turn the entire problem into a the language
of matrix multiplication.  The crucial idea is that checking whether a certain configuration can be excluded is equivalent to checking whether the $\ell^{\infty}$ norm of that configuration interpreted as a vector multiplied by a matrix $C_k\in\mathbb{R}^{(2n)^2\times 4n-k+1}$ exceeds the desired limit.    The $C_k$ indicate all possible $b_i$, $b_j$ pairs that contribute to the interval $I_{a,a+k}\subset (-1/2,1/2)$ for for any $a\in\{-2n,..., 2n-1\}$.  These matrices are precomputed once and stored for a given $n$. 
The only matrix that must be constructed for a bin $b\in B_{n,m}$ is the matrix of refinements of $b$ as $c_b = [c_\gamma]\in \mathbb{R}^{N\times 2n}$, for $N = \prod_{i=1}^n \left( 1 + m b_i \right)$.  We compute the matrix $F\in \mathbb{R}^{N\times (2n)^2}$, where
\begin{eqnarray*}
F[\gamma,n(i-1) + j] = c_\gamma[i]  c_\gamma[j].
\end{eqnarray*}
This means $F C_k \in \mathbb{R}^{N\times 4n-k+1}$ gives the sum of all $c_\gamma*c_\gamma$ along all intervals of length $k$.  A bin $c_\gamma$ is ruled out if 
\begin{eqnarray*}
F[\gamma,\cdot] C_k > 1.28 + \frac{2}{m} + \frac{1}{m^2} 
\end{eqnarray*} 
for any element in the row, for any $k$.  The lower bound can be refined by incorporating the observation from \eqref{eq:refinedLowerBound}.
Obviously, this algorithm is easily parallelizable as each parent bin $b$ can be checked independently of the other parent bins. The actual computation was carried out on \textsc{Omega}, a server at the \textsc{Yale High Performance Computing Center} with 7 CPU nodes each of which had a dedicated GPU. 
The computation took roughly 20.000 CPU hours. \\

\textbf{Acknowledgement.} The first author is supported by NSF Fellowship DMS-1402254, the second author was partially supported by SFB 1060 of the DFG, an AMS-Simons
travel grant and INET Grant \#INO15-00038. The authors are grateful to Andrew Sherman from the \textsc{Yale High Performance Computing Center} for continued support.

\end{document}